\documentclass[twoside,leqno,twocolumn]{article}  
\usepackage{ltexpprt} 
\usepackage{subfig}
\usepackage{float}
\usepackage{amsmath,amsfonts,amssymb}
\usepackage{framed}
\usepackage{graphicx}
\usepackage{hyperref}

\title{\Large Block Tensor Decomposition for Source Apportionment of Air Pollution\thanks{Supported by the NSF DMS-0915100 and the Clarkson's Institute for a Sustainable Environment.}
}
\author{Philip K. Hopke\thanks{Department of Chemical Engineering, Clarkson University, Potsdam, NY 13699.}\\
\and
Maggie Leung\footnotemark[2]\\
\and
Na Li\thanks{Department of Mathematics, Clarkson University, Potsdam, NY 13699.} \\
\and
Carmeliza Navasca\footnotemark[3] \thanks{Corresponding author. Email: cnavasca@clarkson.edu.}\\}

\date{}

\begin{document}

\maketitle


\begin{abstract} \small\baselineskip=9pt 
The ambient particulate chemical composition data with three particle diameter sizes ($2.5$mm $<D_p < 1.15$ mm, $1.15$mm $<D_p< 0.34$mm and $0.34$mm $<D_p <0.1$mm) collected at a major industrial center in Allen Park in Detroit, MI is examined. Standard multiway (tensor) methods like PARAFAC and Tucker tensor decompositions have been applied extensively to many chemical data. However, for multiple particle sizes, the source apportionment analysis calls for a novel multiway factor analysis. We apply the regularized block tensor decomposition to the collected air sample data. In particular, we use the Block Term Decomposition (BTD) in rank-$(L;L;1)$ form to identify nine pollution sources (Fe+Zn, Sulfur with Dust, Road Dust, two types of Metal Works, Road Salt, Local Sulfate, and Homogeneous and Cloud Sulfate).
\end{abstract}

\section{Introduction}
Source apportionment analysis determines the sources of various air pollutants at a particular location.  Through a collected sample of  environmental data, the ambient particulate chemical composition data is acquired and then analyzed. This technique of identifying and quantifying the sources of air pollutants at a receptor location by resolving the mixture of chemical species into the contributions from the individual source types is called receptor modeling when factor analysis used. Factor analysis utilizes matrix methods, like PCA and PMF. In the study of source apportionment of airborne particles, the measured chemical composition data from the collected samples form a matrix in terms of chemical  species and time samples which can be decomposed into two matrices representing sources contributions and sources profiles. Matrix methods have been a very effective tool for source apportionment for collected samples of fine particles as well as coarse particles; see e.g. \cite{Polissar,Hopke3,Hopke4}. However, chemical data can have more attributes, such as particle size distribution which is common in identification problems of air pollution sources \cite{Hopke}. Thus, a two-way analysis is not sufficient. Multiway factor analysis (tensor decompositions) can provide better estimations.

We examine the ambient particulate data \cite{Hopke} collected at Allen Park in Detroit, MI between February and April 2002. It follows the following three-way receptor model:
\begin{eqnarray*} \label{entensor}
x_{ijk}=\sum_{p=1}^P a_{ip}b_{jpk} + e_{ijk}
\end{eqnarray*}
where the measured data $x_{ijk}$ is the concentration value of the $i$th sample of the $k$th particle size range $(0.1-0.34 \mu m, 0.34-1.15 \mu  m, 1.15-2.5 \mu m)$ of the $j$th species, $a_{ip}$ is the $p$th source mass contribution during the time units for the $i$th sample, $b_{ijk}$ is the $j$th species mass fraction of the $k$th particle size range from the $p$th source, $e_{ijk}$ is the measurement error, and $P$ is the total number of independent sources.

In this paper, we reformulate the three-way receptor model in a block tensor structure format and utilize a new regularized tensor algorithm found in \cite{Navasca}. Previous tensor applications in chemometrics or environmental data analysis have been limited to tensor CANDECOMP/PARAFAC (CP) decomposition \cite{JJ,Harshman,Bro,Bro2,Henrion,Henrion1}  or tensor Tucker decomposition \cite{Tucker,Tucker1,Lieven4}. Instead we focus on a more general type of tensor decomposition - Block Term Decomposition (BTD) \cite{Lieven1,Lieven2,Lieven3}, and its application in environmental pollution. In fact, the data model fits the BTD in rank-$(L,L,1)$. 

Here we apply the regularized alternating method (RALS) \cite{Navasca} for the BTD problems.  Then, the method is used to analyze a real-life air sample tensor dataset. We show the efficacy of the regularized alternating algorithm (BTD-RALS) for solving BTD through numerical examples as well as prove a convergence property of the algorithm. In the application part, we show that the air sample data model in \cite{Hopke} follows the BTD format, so that the BTD-RALS method is used to identify all the different factor pollution sources.

The paper is organized as follows. We introduce the notation in Section 2. In Sections 3 and 4, we describe some background of the block term decomposition and then discuss a convergence property of the algorithm. The real data example is given in Section 6. Concluding remarks are discussed in Section 6.

\section{Preliminaries} 
We denote the scalars in $\mathbb{R}$ with lower-case letters $(a,b,\ldots)$ and the vectors with bold lower-case letters $(\bf{a},\bf{b},\ldots)$.  The matrices are written as bold upper-case letters $(\bf{A}, \bf{B},\ldots)$ and the symbol for tensors are calligraphic letters $(\mathcal{A},\mathcal{B},\ldots)$. The subscripts represent the following scalars:  $\mathcal{(A)}_{ijk}=a_{ijk}$, $(\bold{A})_{ij}=a_{ij}$, $(\bold{a})_i=a_i$ and $\bold{A}_r$ is the $r$th column of $\bold{A}$. The superscripts indicate the length of the vector or the size of the matrices. For example, $\bold{b}^{K}$ is a vector with length $K$ and $\bold{B}^{N \times K}$ is a $N \times K$ matrix. In addition, the lower-case superscripts on a matrix indicate the mode in which it has been matricized.  For example, $\bold{T}_{(n)}$ is the mode-$n$ matricization of the tensor $\mathcal{T} \in \mathbb{R}^{I \times J \times K}$ for $n=1,2,3$.

\begin{Definition}
The Kronecker product of matrices $\bold{A}$ and $\bold{B}$ is defined as
{\small \begin{eqnarray*}
\bold{A} \otimes \bold{B}= \left[
\begin{array}{ccc}
a_{11}\bold{B}  & a_{12}\bold{B} & \ldots \\
a_{21}\bold{B}  & a_{22}\bold{B} & \ldots \\
\vdots                  & \vdots     & \ddots
\end{array} \right].
\end{eqnarray*}}
\end{Definition}

\begin{Definition}
The Khatri-Rao product  is the ``matching columnwise" Kronecker product. Given matrices $\bold{A} \in \mathbb{R}^{\mathnormal{I} \times \mathnormal{K}}$ and $\bold{B} \in \mathbb{R}^{\mathnormal{J} \times \mathnormal{K}}$, their Khatri-Rao product is denoted by $\bold{A} \odot \bold{B}$. The result is a matrix of size $(\mathnormal{IJ} \times \mathnormal{K})$ defined by
\begin{eqnarray*}
\bold{A} \odot \bold{B}= [\bold{A}_1 \otimes \bold{B}_1~~ \ldots~~ \bold{A}_K \otimes \bold{B}_K].
\end{eqnarray*}
If $\bold{a}$ and $\bold{b}$ are vectors, then the Khatri-Rao and Kronecker products are identical, i.e., $\bold{a} \otimes \bold{b} = \bold{a} \odot \bold{b}$.
\end{Definition}

\begin{Definition}
Let $\bold{A} = [\bold{A}_1~\ldots~\bold{A}_R]$ and $\bold{B}=[\bold{B}_1~\ldots~\bold{B}_R]$ be two partitioned matrices. Then we define a product of $\bold{A}$ and $\bold{B}$, denoted $\odot_p$, which is
$$\bold{A} \odot_p \bold{B} = [\bold{A}_1 \otimes \bold{B}_1~\ldots~\bold{A}_R \otimes \bold{B}_R].$$
\end{Definition}

\begin{Definition}[Mode-$n$ matricization]
Matricization is the process of reordering the elements of an $N$th order tensor into a matrix. The mode-$n$ matricization of a tensor $\mathcal{T} \in \mathbb{R}^{I_1 \times I_2 \times \cdots \times I_N}$ is denoted by $\bold{T_{(n)}}$ and arranges the mode-$n$ fibers, the vectors obtained from fixing every index with the exception of the $n$th mode,  as the columns of the resulting matrix.
\end{Definition}

\begin{Definition}[Vectorization]
The vectorization of a matrix $\bold{M}=[\bold{m}_1~\cdots ~ \bold{m}_n] \in \mathbb{R}^{m \times n}$, where $\bold{m}_i$ is the $i$th column of $\bold{M}$, is denoted by $vec(\bold{M})$ which is a vector of size $mn$ defined by
\begin{eqnarray*}
vec(\bold{M}) =\left[ \begin{array}{c} \bold{m}_1 \\ \vdots \\ \bold{m}_n \end{array} \right] .
\end{eqnarray*}
\end{Definition}
 
\begin{Definition}[Frobenius-norm]
The Frobenius norm of a tensor $\mathcal{X} \in \mathbb{R}^{I_1 \times I_2 \times \cdots \times I_N}$ is the square root of the sum of the squares of all its elements. The formula is
\begin{eqnarray*}
\Vert \mathcal{X} \Vert_F = \sqrt{\sum_{i_1=1}^{I_1} \sum_{i_2=1}^{I_2}\cdots \sum_{i_N=1}^{I_N} x_{i_1i_2\cdots i_N}^2}.
\end{eqnarray*}
\end{Definition}

\section{Block Term Decompositions}
Let $\mathcal{X}$ be a real-valued third-order tensor of size $I \times J \times K$. A decomposition of $\mathcal{X}$ in a sum of rank-$(L_r,L_r,1)$ terms with $1\leq r \leq R$ is a decomposition of the form
\begin{eqnarray}\label{eq:btd1}
\mathcal{X} = \sum_{r=1}^R \bold{E}_r \circ \bold{c}_r,
\end{eqnarray}
in which the rank of the matrices $\bold{E}_r \in \mathbb{R}^{ I \times J}$ is $L_r$ and $\bold{c}_r $ are vectors of length K. Elementwise, the decomposition is defined as 
$$x_{ijk}=\sum_{r=1}^R e^{(r)}_{ij}\cdot c^{(r)}_k.$$

So $\mathcal{X}$ is decomposed into a sum of matrix-vector outer products. If we decomposes $\bold{E}_r$ into two matrices, i.e., $\bold{E}_r = \bold{A}_r \cdot \bold{B}_r^{{T}}$, where the matrix $\bold{A}_r \in \mathbb{R}^{I \times L_r}$ and the matrix $\bold{B}_r \in \mathbb{R}^{J \times L_r}$ are rank $L_r$, then the equation (\ref{eq:btd1}) can be written as the following formula,
\begin{eqnarray}\label{eq:btd2}
\mathcal{X} = \sum_{r=1}^R (\bold{A}_r \cdot \bold{B}_r^{T} ) \circ \bold{c}_r . 
\end{eqnarray}

We define $\bold{A} = [\bold{A}_1~\cdots~ \bold{A}_R]$, $\bold{B}=[\bold{B}_1~\cdots~\bold{B}]$, $\bold{C} = [\bold{c}_1~\cdots~\bold{c}_R]$ and call them the factor matrices of the BTD-$(L_r,L_r,1)$, where $\bold{A}_r$ and $\bold{B}_r$ are submatrices with $r=1,\dots,R$ . When $L_r \equiv L$ for $1 \leq r \leq R$, the decomposition is called BTD in rank-$(L,L,1)$. We will focus on this type of factorization and analyze the environmental dataset with BTD-$(L,L,1)$.

\begin{figure}[h!]
\includegraphics[width=.5\textwidth]{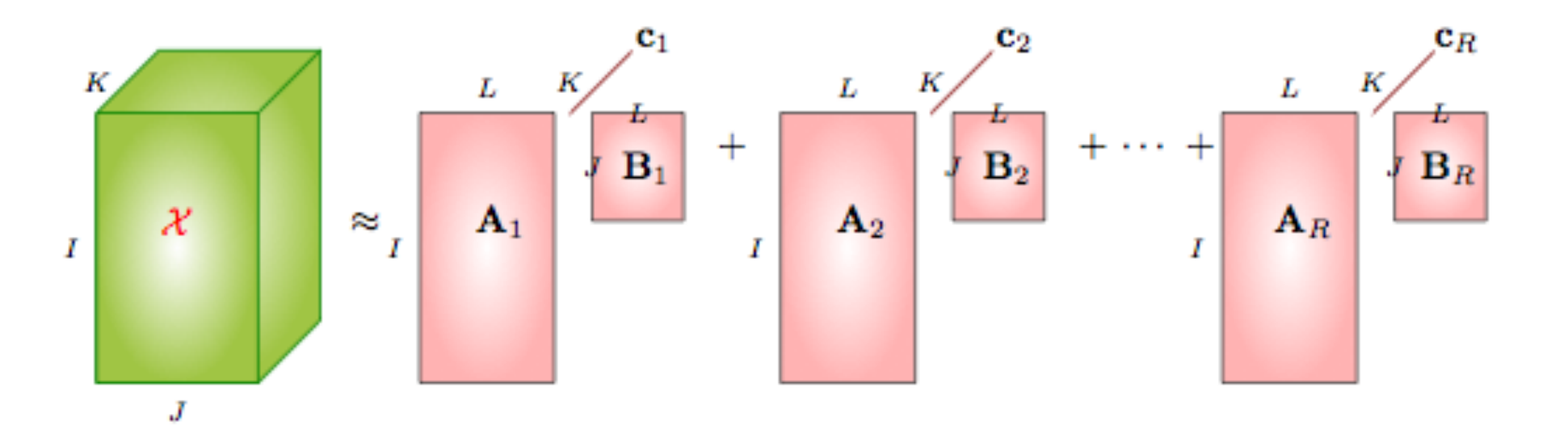}
\caption{\bf BTD-$(L,L,1)$ for $\mathcal{X} \in \mathbb{R}^{I \times J \times K}$}
\label{fig:btd}
\end{figure}

See Figure \ref{fig:btd} for the BTD in rank-$(L,L,1)$ for a third-order tensor.

So, for a given tensor $\mathcal{X} \in \mathbb{R}^{I \times J \times K}$, the BTD-$(L,L,1)$ will minimize the error function
\begin{eqnarray}\label{eq:cost0}
\qquad f(\bold{A},\bold{B},\bold{C})=\biggl\Vert \mathcal{X}- \sum_{r=1}^R(\bold{A}_r \cdot \bold{B}_r^T)\circ \bold{c}_r\biggl\Vert_F^2.
\end{eqnarray}

The general case of the block term decomposition is called BTD in rank-$(L,M,N)$. It decomposes a tensor $\mathcal{X} \in \mathbb{R}^{I \times J \times K}$ in a sum of rank-$(L,M,N)$ terms of the form
\begin{eqnarray}\label{eq:btdlmn}
\mathcal{X} = \sum_{r=1}^R \mathcal{D}_r \times_1 \bold{A}_r \times_2 \bold{B}_r \times_3 \bold{C}_r,
\end{eqnarray}
in which $\mathcal{D}_r \in \mathbb{R}^{L\times M \times N}$ and where $\bold{A}_r \in \mathbb{R}^{I \times L}$, $\bold{B}_r \in \mathbb{R}^{J \times M}$ and $\bold{C}_r \in \mathbb{R}^{K \times N}$ are full column rank, $1 \leq r \leq R$.

\begin{figure}[h!]
\includegraphics[width=.5\textwidth]{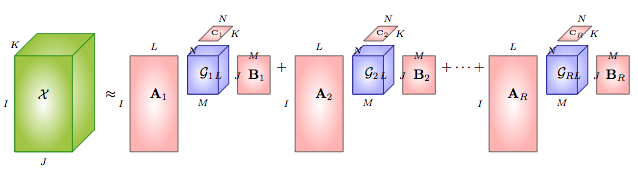}
\caption{\bf BTD-$(L,M,N)$ for $\mathcal{X} \in \mathbb{R}^{I \times J \times K}$}
\label{fig:btdlmn}
\end{figure}

The product `$\times_1$', `$\times_2$' and `$\times_3$' are called Tucker product which defines the multiplication between a tensor and a matrix. Figure \ref{fig:btdlmn} shows the BTD in rank-$(L,M,N)$ for a third-order tensor.

Actually, each term of the equation (\ref{eq:btdlmn}) is a Tucker model \cite{Tucker}. We see here if let $N=1$, then each core tensor becomes a matrix, thereby reducing the general model to BTD-$(L_r, L_r, 1)$. If the ranks of the matrices for each term are same, then it becomes the BTD-$(L,L,1)$. Furthermore, if $L=M=N=1$, then the each core tensor is just a scalar and all the matrices $\bold{A}_r$, $\bold{B}_r$ and $\bold{C}_r$ are vectors. So, the BTD-$(1,1,1)$ is exactly a CP decomposition,



\begin{eqnarray}\label{eq:cp}
\mathcal{X}=\sum_{r=1}^R \bold{a}_r \circ \bold{b}_r \circ \bold{c}_r.
\end{eqnarray}
The BTD in rank-$(L,L,1)$ is a general case of the CP decomposition. Thus, the standard and regularized algorithms for solving for CP can be applied to the BTD case. More details on the algorithms are discussed in the next section.

\section{Algorithm for Block Term Decomposition}

For the BTD-$(L,L,1)$, equation (\ref{eq:btd2}) can be expressed to three equivalent equations. If we take the three different modes matricization on equation (\ref{eq:btd2}), then we have 
\begin{eqnarray*}
\bold{X}_{(1)} &=& \bold{A} (\bold{C} \odot_p \bold{B})^{T},  \\
\bold{X}_{(2)} &=& \bold{B} (\bold{C} \odot_p \bold{A})^{T},  \\
\bold{X}_{(3)} &=& \bold{C} [(\bold{B}_1 \odot \bold{A}_1) \bold{1}_L ~ \cdots ~ (\bold{B}_R \odot \bold{A}_R)\bold{1}_L]^{T},
\end{eqnarray*}
where $\bold{1}_L$ is a column vector of all ones of length $L$.

The paper of De Lathauwer and Nion \cite{Lieven3} proposes an algorithm to solve the block term decomposition in rank $(L,L,1)$ for a third-order tensor, called the alternating least-squares for BTD (BTD-ALS). 

Given a third-order tensor $\mathcal{X} \in \mathbb{R}^{I \times J \times K}$, our problem is to minimize the function (\ref{eq:cost0}), i.e., 

\begin{eqnarray}\label{eq:problem}
\min_{\bold{A}, \bold{B}, \bold{C}} \biggl\Vert \mathcal{X} - \sum_{r=1}^R (\bold{A}_r \cdot \bold{B}_r^T) \circ \bold{c}_r\biggl\Vert_F^2,
\end{eqnarray}
with respect to the factor matrices $\bold{A}$, $\bold{B}$ and $\bold{C}$.

From the three equations above, we obtain the following three expressions for (\ref{eq:problem}):
\begin{eqnarray*}
&&\min_{\bold{A}, \bold{B}, \bold{C}} \Vert \bold{X}_{(1)} - \bold{A}(\bold{C} \odot_p \bold{B})^T\Vert_F^2, \\
&& \min_{\bold{A}, \bold{B}, \bold{C}} \Vert \bold{X}_{(2)} - \bold{B} (\bold{C} \odot_p \bold{A})^T\Vert_F^2, \\
&& \mbox{and}\\
&& \min_{\bold{A}, \bold{B}, \bold{C}} \Vert \bold{X}_{(3)} - \bold{C}[(\bold{B}_1 \odot \bold{A}_1) \bold{1}_L \cdots (\bold{B}_R \odot \bold{A}_R)\bold{1}_L]^T\Vert_F^2.
\end{eqnarray*}

These three are equivalent. Instead of solving (\ref{eq:problem}) for the three factors one time, we can use these three equations by fixing all factor matrices but one factor each time. Then the problem reduces to three coupled linear least-squares subproblems. We have
\begin{eqnarray*}
&\bold{A}^{k+1}&=\displaystyle\mathop{\mathrm{argmin}}_{\widehat{\bold{A}}\in \mathbb{R}^{I\times LR}}\Vert\bold{X_{(1)}}-\widehat{\bold{A}}  (\bold{C}^{k}\odot_p \bold{B}^{k})^T\Vert_F^2 , \\
&\bold{B}^{k+1}&=\displaystyle\mathop{\mathrm{argmin}}_{\widehat{\bold{B}}\in \mathbb{R}^{J\times LR}}\Vert\bold{X_{(2)}}-\widehat{\bold{B}}(\bold{C}^{k}\odot_p \bold{A}^{k+1})^T\Vert_F^2,\\
&\bold{C}^{k+1}&=\displaystyle\mathop{\mathrm{argmin}}_{\widehat{\bold{C}}\in \mathbb{R}^{K\times R}}\Vert\bold{X_{(3)}}- \widehat{\bold{C}}[(\bold{B}^{k+1}_1 \odot \bold{A}^{k+1}_1) \bold{1}_L \\
&& \qquad\qquad\qquad \cdots (\bold{B}^{k+1}_R \odot \bold{A}^{k+1}_R) \bold{1}_L]^T\Vert_F^2 , 
\end{eqnarray*}
where $\bold{X}_{(1)} \in\mathbb{R}^{I \times JK}$, $\bold{X}_{(2)} \in \mathbb{R}^{J \times IK}$ and $\bold{X}_{(3)} \in \mathbb{R}^{K \times IJ}$ are the mode-1, mode-2 and mode-3 matricizations of tensor $\mathcal{X}$.

Given the initials $\bold{A}^0$, $\bold{B}^0$ and $\bold{C}^0$, then at the $(k+1)th$ iteration, we hold $\bold{B}^k$ and $\bold{C}^k$ to update the factor $\bold{A}$ to get $\bold{A}^{k+1}$, then $\bold{A}^{k+1}$ and $\bold{C}^k$ are held to update $\bold{B}$ and obtain $\bold{B}^{k+1}$. Similarly, we hold $\bold{A}^{k+1}$ and $\bold{B}^{k+1}$ to obtain $\bold{C}^{k+1}$. Usually, the Frobenius norm of the error between the given tensor and the updated tensor is measured at each iteration to provide a stopping criterion. 

There are some disadvantages to alternating least-squares algorithm (see \cite{Kolda}, \cite{Na}). This method is not guaranteed to converge to a global minimum or even a stationary point of the cost function (\ref{eq:problem}), only to a solution where the objective function ceases to decrease. Another issue of this method is that sometimes it needs a high number of iterations to converge, creating a swamp. In order to remove the swamp, \cite{Navasca} introduced a regularization method with the regularization parameter $\lambda_n$. The new algorithm for BTD-$(L,L,1)$ is called BTD-RALS method.

We add the regularization item for each subproblem in the above three equations,
\begin{eqnarray*}
\bold{A}^{n+1} &=& \displaystyle\mathop{\mathrm{argmin}}_{\widehat{\bold{A}}\in \mathbb{R}^{I\times LR}}\Vert\bold{X_{(1)}}-\widehat{\bold{A}}  (\bold{C}^{n}\odot_p\bold{B}^{n})^T\Vert_F^2 \\
&& + \lambda_n \Vert \widehat{\bold{A}}-\bold{A}^n\Vert_F^2, \\
\bold{B}^{k+1}&=&\displaystyle\mathop{\mathrm{argmin}}_{\widehat{\bold{B}}\in \mathbb{R}^{J\times LR}}\Vert\bold{X_{(2)}}-\widehat{\bold{B}}(\bold{C}^{n}\odot_p\bold{A}^{n+1})^T\Vert_F^2 \\
&&+ \lambda_n \Vert \widehat{\bold{B}}-\bold{B}^n\Vert_F^2,\\
\bold{C}^{k+1}&=&\displaystyle\mathop{\mathrm{argmin}}_{\widehat{\bold{C}}\in \mathbb{R}^{K\times R}}\Vert\bold{X_{(3)}}-\widehat{\bold{C}}[(\bold{B_1}^{n+1}\odot\bold{A_1}^{n+1})\bold{1}_L, \dots,  \\
&&(\bold{B_R}^{n+1} \odot \bold{A_R}^{n+1})\bold{1}_L]^T\Vert_F^2 + \lambda_n \Vert\widehat{\bold{C}} - \bold{C}^n\Vert_F^2. 
\end{eqnarray*}

In alternating fashion, these three subproblems are solved for the factor matrices $\bold{A}$, $\bold{B}$ and $\bold{C}$. The regularization parameters $\lambda_n$, $n=0, 1, \dots$, are given by a nonnegative decreasing sequence and at each iteration the parameters are the same for the three updated factor matrices. The rules for choosing the regularization parameters is also discussed in \cite{Navasca}. The regularized algorithm is summarized in the following Table \ref{tab:btdll1rals}. The number of iterations $N$ in the algorithm is set to a large number, and a stopping criterion can be used.
\begin{table}[h!]
\begin{framed}
  \textbf{\textsf{RBTD-$(L,L,1)$-Algorithm}} \cite{Navasca}
\vskip5pt  
\textbf{procedure} RBTD-$(L,L,1)$($\mathcal{X}, R, N, \lambda_n$)
\vskip5pt  
\hskip10pt give initial guess $\bold{A}^0 \in \mathbb{R}^{I \times R}$, $\bold{B}^0 \in \mathbb{R}^{J \times R}$, $\bold{C}^0\in \mathbb{R}^{K \times R}$, $\lambda_0$
\vskip5pt  
\hskip10pt \textbf{for} $n=1,\dots, N$ \textbf{do}
\vskip5pt  
\hskip26pt  $\bold{W} \leftarrow [(\bold{C}^n \odot_p \bold{B}^n); \lambda_n\bold{I}^{LR \times LR}] \in \mathbb{R}^{(JK+LR) \times LR}$
\vskip5pt  
\hskip26pt  $\bold{S} \leftarrow [\bold{X}_{(1)} ; \lambda_n(\bold{A}^n)^{T}] \in \mathbb{R}^{(JK+LR) \times I}$
\vskip5pt  
\hskip26pt $\bold{A}^{n+1} \leftarrow \bold{W}/ \bold{S}$ ------ \% solving least squares to update $\bold{A}$
\vskip10pt

\hskip26pt  $\bold{W} \leftarrow [(\bold{C}^n \odot_p \bold{A}^{n+1}); \lambda_n\bold{I}^{LR \times LR}] \in \mathbb{R}^{(IK+LR) \times LR}$
\vskip5pt  
\hskip26pt  $\bold{S} \leftarrow [\bold{X}_{(2)} ; \lambda_n(\bold{B}^n)^{T}] \in \mathbb{R}^{(IK+LR) \times J}$
\vskip5pt  
\hskip26pt $\bold{A}^{n+1} \leftarrow \bold{W}/ \bold{S}$ ------ \% solving least squares to update $\bold{B}$
\vskip10pt

\hskip26pt $\bold{W} \leftarrow [((\bold{B}_1^{n+1} \odot \bold{A}_1^{n+1})\bold{1}_L, \dots, (\bold{B}_R^{n+1} \odot \bold{A}_R^{n+1})\bold{1}_L); \lambda_n\bold{I}^{R\times R}] \in \mathbb{R}^{(IJ+R) \times R}$
\vskip5pt  
\hskip26pt  $\bold{S} \leftarrow [\bold{X}_{(3)}; \lambda_n(\bold{C}^n)^{T}] \in \mathbb{R}^{(IJ+R)\times K} $
\vskip5pt  
\hskip26pt $\bold{C}^{n+1} \leftarrow \bold{W}/\bold{S}$ ------ \% solving least squares to update $\bold{C}$
\vskip10pt

\hskip10pt \textbf{end for}
\vskip5pt  
\hskip10pt return $\bold{A}^{N}$, $\bold{B}^{N}$, $\bold{C}^{N}$
\vskip5pt  
\textbf{end procedure}
\end{framed}
\caption{\bf Regularized algorithm of BTD-$(L,L,1)$ with rank $R$ for a third-order tensor $\mathcal{X} \in \mathbb{R}^{I \times J \times K}$}
\label{tab:btdll1rals}
\end{table}

Numerically, such regularized method is more efficient than the BTD-ALS in terms of reducing the number of iterations and accelerating the convergence.

\begin{Example}[Numerical example with a swamp]
Here we give an example to show the swamp reducing technique of the BTD-RALS. 

We create a tensor $\mathcal{X} \in \mathbb{R}^{10 \times 15 \times 28}$, which satisfies a block term decomposition in rank-$(3,3,1)$ with $R=3$. The factor matrices are $\bold{A} \in \mathbb{R}^{10 \times 9}$, $\bold{B} \in \mathbb{R}^{15 \times 9}$ and $\bold{C} \in \mathbb{R}^{28 \times 3}$, and the factorization equation is $$\bold{X} =\sum_{r=1}^{3} (\bold{A}_r \cdot \bold{B}_r^{\text{T}})\circ \bold{c}_r.$$

We use the same random initials for both BTD-ALS and BTD-RALS methods. Figure \ref{fig:btdnumerical2} shows that the BTD-RALS algorithm only takes 1558 iterations to reach an error within $10^{-4}$, however, the BTD-ALS algorithm does not decrease the error after 20,000 iterations.

\begin{figure}[h!]
\begin{center}
\includegraphics[width=.5\textwidth]{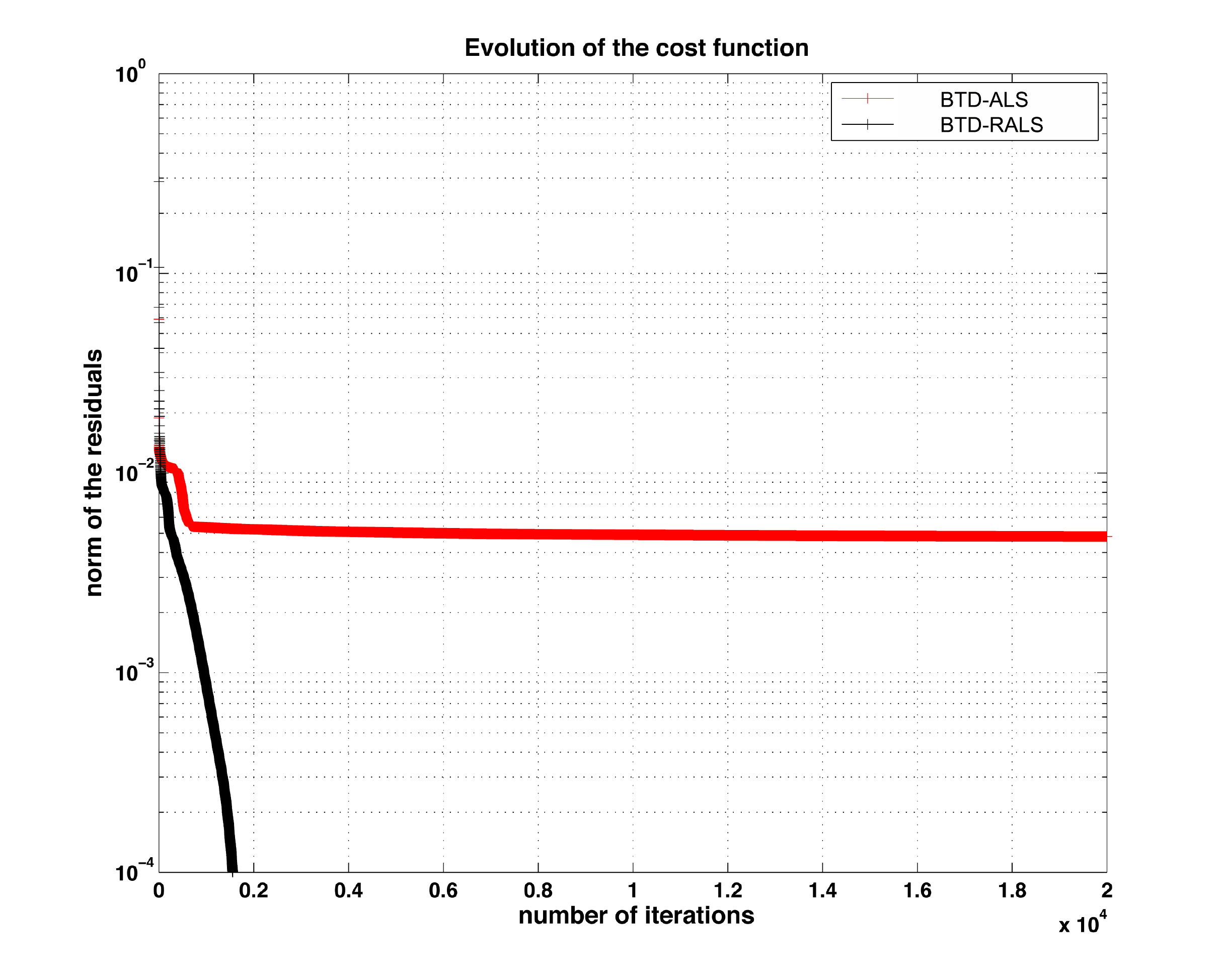}
\caption{\bf The comparison of the BTD and RBTD with the same initials}
\label{fig:btdnumerical2}
\end{center}
\end{figure}
\end{Example}

\subsection{Convergence Property of RALS}
The regularized method of ALS (RALS) \cite{Navasca} for solving tensor CP decomposition and the convergence property of RALS have been studied in \cite{Na}. We have pointed out that the decomposition BTD-$(L,L,1)$ is a general case of CP decomposition. In this section, we will show that the BTD-RALS has the same framework with the RALS and thus, the BTD-RALS has the same convergence property. 

We can view the BTD-$(L,L,1)$ problem as a nonlinear optimization. Problem (\ref{eq:problem}) has the following expression,
\begin{eqnarray}\label{eq:problem1}
&&\min_{\bold{A}, \bold{B}, \bold{C}} f(\bold{A}, \bold{B}, \bold{C}) \\
&=&\min_{\bold{A}, \bold{B}, \bold{C}} \sum_{i=1}^I \sum_{j=1}^J \sum_{k=1}^K \biggl[x_{ijk} - \sum_{r=1}^R\biggl(\sum_{l=1}^L a_{il}^{(r)}b_{jl}^{(r)} \biggl)c_k^{(r)}\biggl]^2, \notag
\end{eqnarray}
where $\bold{A}=[\bold{A}_1\cdots \bold{A}_R]$, $\bold{B} = [\bold{B}_1 \cdots \bold{B}_R]$ and $\bold{C} = [\bold{c}_1 \cdots \bold{c}_R]$ are the factor matrices. $a_{il}^{(r)}$ denotes the $il$ element ($i$th row and $l$th column) of the matrix $\bold{A}_r$, $b_{jl}^{(r)}$ expresses the $jl$ element of the matrix $\bold{B}_r$, and $c_k^{(r)}$ is the $k$th element in the vector $\bold{c}_r$.

So the cost function can be seen as a function from $\bold{y}$ to $\mathbb{R}$, where
$$\bold{y}=vec([vec(\bold{A})~vec(\bold{B})~vec(\bold{C})]) \in \mathbb{R}^n,$$
with $n=(IL+JL+K)R$.

Let $vec(\bold{A}) = \bold{y}_1 \in \mathbb{R}^{ILR}$, $vec(\bold{B}) = \bold{y}_2 \in \mathbb{R}^{JLR}$ and $vec(\bold{C}) = \bold{y}_3 \in \mathbb{R}^{KR}$. We partition the vector $\bold{y} \in \mathbb{R}^n$ into 3 component vectors $\bold{y}_i \in \mathbb{R}^{n_i}$, $i=1,2,3$, where $n_1 = ILR$, $n_2 = JLR$ and $n_3 = KR$. It follows that $\bold{y}=\bold{y}_1 \times \bold{y}_2 \times \bold{y}_3 \in \mathbb{R}^{n_1} \times \mathbb{R}^{n_2} \times \mathbb{R}^{n_3} = \mathbb{R}^n$. Thus, the BTD-$(L,L,1)$ can be reformulated to the following problem,
\begin{eqnarray*}
\mbox{\bf minimize}&&\quad f(\bold{y}_1,\bold{y}_2,\bold{y}_3), \\
\mbox{\bf subject to}&& \quad \bold{y} \in  \mathbb{R}^{n_1} \times \mathbb{R}^{n_2} \times \mathbb{R}^{n_3} =\mathbb{R}^n. 
\end{eqnarray*}

Therefore, the BTD-RALS method for solving the BTD-$(L,L,1)$ updates each component of $\bold{y}$ in turn. Starting from a given initial point $\bold{y}^0 = vec([vec(\bold{A}^0)~vec(\bold{B}^0)~vec(\bold{C}^0)]) \in \mathbb{R}^n$, a sequence $\{ ( \bold{y}_1^k, \bold{y}_2^k, \bold{y}_3^k ) \}$ is generated by the following equations
\begin{eqnarray*}
\bold{y}_1^{k+1} &=& \displaystyle\mathop{\mathrm{argmin}}_{\bold{z} \in \mathbb{R}^{n_1}} \{ f(\bold{z},  \bold{y}_2^{k}, \bold{y}_3^k) + \lambda_k \Vert \bold{y}_1^k- \bold{z}\Vert^2 \} ,  \\
\bold{y}_2^{k+1} &=& \displaystyle\mathop{\mathrm{argmin}}_{\bold{z} \in \mathbb{R}^{n_2}}\{ f(\bold{y}_1^{k+1}, \bold{z}, \bold{y}_3^k) + \lambda_k \Vert \bold{y}_2^k- \bold{z}\Vert^2, \\
\bold{y}_3^{k+1} &=& \displaystyle\mathop{\mathrm{argmin}}_{\bold{z} \in \mathbb{R}^{n_3}}\{ f(\bold{y}_1^{k+1}, \bold{y}_2^{k+1} , \bold{z}) + \lambda_k \Vert \bold{y}_3^k - \bold{z} \Vert^2 \}. 
\end{eqnarray*}

Thus, the BTD-RALS algorithm has the same framework as the RALS for CP decomposition \cite{Na}. Moreover, the convergence property of RALS is applied to the BTD-RALS. 

For the BTD-RALS method, we have the following theorem,
\begin{theorem}
Suppose that the sequence $\{(\bold{A}^k, \bold{B}^k, \bold{C}^k)\}$ obtained from the BTD-RALS has limit points. Then every points $(\overline{\bold{A}}, \overline{\bold{B}}, \overline{\bold{C}})$ is a critical point of the problem (\ref{eq:problem1}).
\end{theorem}

Recall that a sequence that has a limit point $(\overline{\bold{A}}, \overline{\bold{B}}, \overline{\bold{C}})$ means that there exists a subsequence of $\{(\bold{A}^k, \bold{B}^k, \bold{C}^k)\}$ that converges to $\{(\bold{A}^k, \bold{B}^k, \bold{C}^k)\}$.  The following is a critical point definition found in \cite{Bertsekas}. The critical point $(\overline{\bold{A}}, \overline{\bold{B}}, \overline{\bold{C}})$ of the problem (\ref{eq:problem1}) is a point such that $$\bigtriangledown f(\overline{\bold{A}}, \overline{\bold{B}}, \overline{\bold{C}})^T\biggl((\bold{A}, \bold{B}, \bold{C}) - (\overline{\bold{A}}, \overline{\bold{B}}, \overline{\bold{C}})\biggl) \geq 0, \; \forall (\bold{A}, \bold{B}, \bold{C}).$$

According to this theorem, we can see that for a non-degenerate BTD problem \cite{Kolda}, \cite{Na}, the limit points obtained from the BTD-RALS method are the critical point of the original cost function (\ref{eq:cost0}). 


\section{Experiments}
In this section, the air pollution collected data \cite{Hopke} is analyzed in rank-$(9, 9, 1)$ BTD using the BTD-RALS algorithm.  Several figures are then created from the numerical results which explain the sources' identification via the sources profiles and time series of the source contributions. We also provide some numerical simulations on randomly generated tensor noisy data.

%

\subsection{Environmental Data}
The original data was collected from 25 February to 10 April 2002 at Allen Park in Detroit, Michigan. The data was then sampled by using a three-stage Davis Rotating-drum Universal-size-cut Monitoring (DRUM) impactor \cite{Raabe}, \cite{Hopke}. The particles were collected in three size modes, $2.5\mu m > D_p > 1.15\mu m$, $1.15\mu m > D_P > 0.34\mu m$ and $0.34 \mu m > D_P > 0.1 \mu m$, where $D_p$ denotes the particle diameter. The air sample was also analyzed for elements of higher atomic number. The 27 chemical species found were Na, Mg, Al, Si, P, S, Cl, K, Ca, Ti, V, Cr, Mn, Fe, Co, Ni, Cu, Zn, Ga, As, Se, Br, Rb, Sr, Zr, Mo, and Pb.  

The data we study is considered as a function of size, time and chemical composition (a.k.a. elemental species). If we use $i$ to denote chemical species, $j$ to express particle size and $k$ to be the time sample, then a data point $x_{ijk}$ can be expressed as the concentration value of the $i$th chemical species of the $j$th particle size of the $k$th time sample.

\begin{figure}[h!]
\begin{center}
\includegraphics[width=0.53\textwidth]{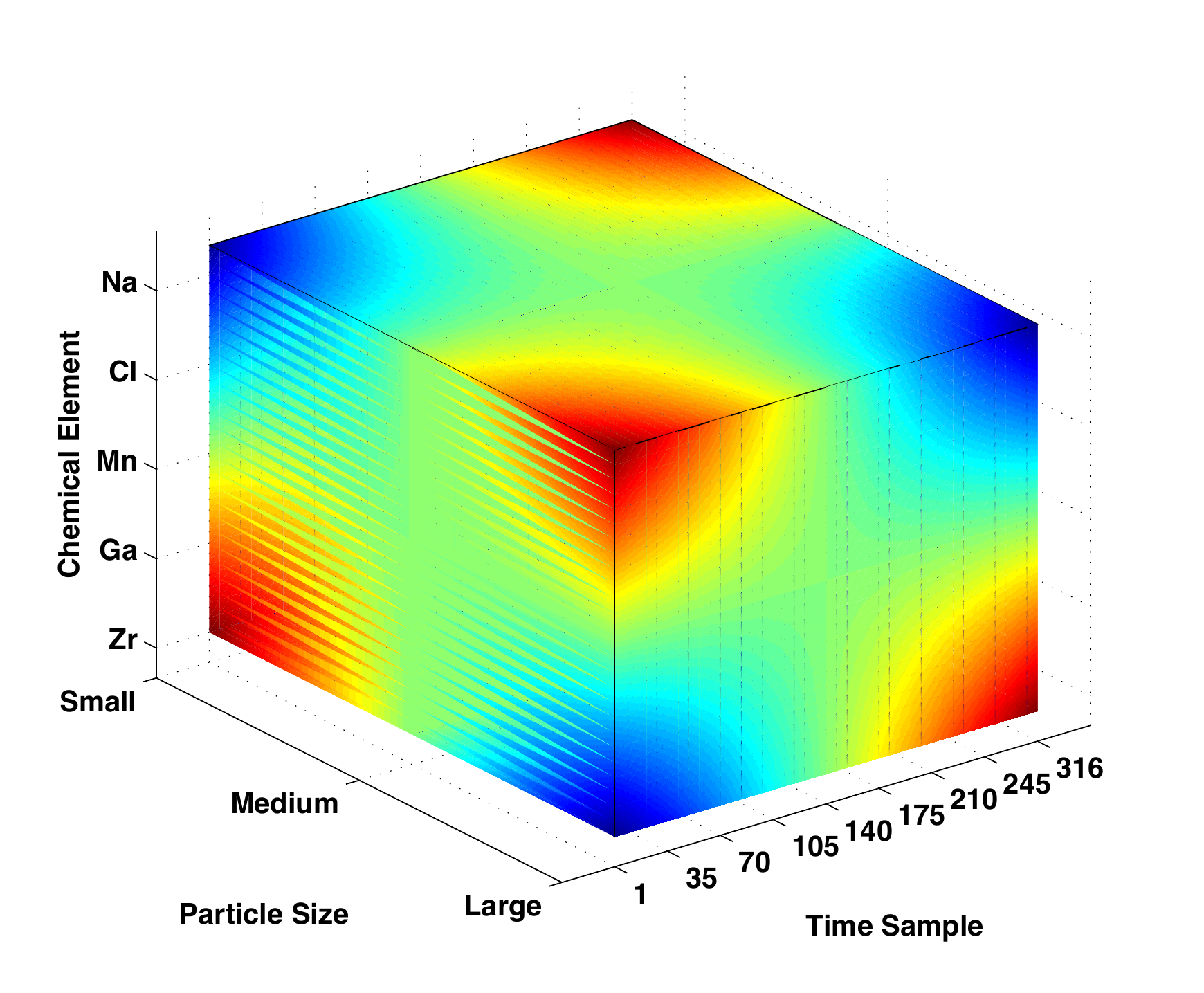}
\end{center}
\caption{\bf The sampled data are placed on the three planes to construct the air tensor $\mathcal{X}$. Particle size denotes mode-$J$, Chemical composition denotes mode-$I$ and Time sample denotes mode-$K$. So, the element $x_{ijk}$ in $\mathcal{X}$ is the concentration value of the $i$th chemical species of the $j$th particle size range of the $k$th time sample.}
\label{fig:figure1}
\end{figure}

See Figure \ref{fig:figure1} for the air sample tensor picture.

There are $27$ chemical species and three different size of particle (small, medium, large), and according to \cite{Hopke}, the time sample are $316$ (3 time samples at first day; 8 time points on each day from the second to the last second day; 1 time point in the last day). So, the air tensor is a third-order tensor of size $27\times 3\times 316$.

\subsection{Model Description} 
According to \cite{Hopke}, in order to separate the different factor sources from the dataset, the main equation is as follows:
\begin{eqnarray}\label{eq:model1}
x_{ijk}=\sum_{p=1}^P a_{ijp}b_{kp} + e_{ijk},
\end{eqnarray}
where the $x_{ijk}$ is the element of the third-order tensor obtained from the previous section, i.e., the concentration value of the $i$th chemical species of the $j$th particle size range of the $k$th time sample. In the elemental form, $b_{kp}$ is the $p$th source mass contribution during the time units for the $k$th sample, $a_{ijp}$ is the $i$th species mass fraction of the $j$th particle size range from the $p$th source. $e_{ijk}$ is the residual associated with the $i$th species concentration measured in the $k$th sample of the $j$th size range, and $P$ is the total number of independent sources. 

Since each entry $x_{ijk}$ denotes the concentration, $a_{ijp}$ is the species mass fraction and $b_{kp}$ is the source contribution, then we need add non-negativity constraints on $a_{ijp}$ and $b_{kp}$ when the tensor decomposition is applied. 

We see that in the equation (\ref{eq:model1}), $x_{ijk}$ is an element of the air sample tensor $\mathcal{X}$. For each fixed $p$, $a_{ijp}$ is an entry of a matrix $\bold{A}_p$ and $b_{kp}$ is an element of a vector $\bold{b}_{p}$. Finally, the error tensor can be denoted as $\mathcal{E}$. Therefore, the equation (\ref{eq:model1}) is equivalent to the following form,
\begin{eqnarray}\label{eq:modeln}
\mathcal{X} = \sum_{p=1}^P \bold{A}_p \circ \bold{b}_p + \mathcal{E}.
\end{eqnarray}

Since $\mathcal{X} \in \mathbb{R}^{27 \times 3 \times 316}$, for each matrix $\bold{A}_p \in \mathbb{R}^{27 \times 3}$, it can be decomposed into two matrices $\bold{C}_p \in \mathbb{R}^{27 \times L_p}$ and $\bold{B} \in \mathbb{R}^{3 \times L_p}$ so that $\bold{A}_p = \bold{C}_p \cdot \bold{D}_p^T$, where $L_p$ is the rank of matrix $\bold{A}_p$. Therefore, the model also can be written as follows,
\begin{eqnarray}\label{eq:modeln1}
\mathcal{X} = \sum_{p=1}^P(\bold{C}_p \cdot \bold{D}_p^T) \circ \bold{b}_p + \mathcal{E}.
\end{eqnarray}

Thus, our goal is to find the matrices $\bold{C} = [\bold{C}_1~\cdots ~ \bold{C}_P]$, $\bold{D}=[\bold{D}_1~\cdots ~\bold{D}_P]$ and $\bold{B}=[\bold{b}_1 ~ \cdots ~ \bold{b}_P]$ to minimize the following function, 
\begin{eqnarray}\label{eq:costfunction}
\qquad\quad Q(\bold{C}, \bold{D}, \bold{B}) = \biggl \Vert \mathcal{X} - \sum_{p=1}^P(\bold{C}_p \cdot \bold{D}_p^T)\circ \bold{b}_p \biggl \Vert_F^2.
\end{eqnarray}
This is exactly the error function for BTD (\ref{eq:cost0}) in rank-$(L_r, L_r, 1)$. When $L_i = L_j$, $1 \leq i, j \leq R$, it is the error function for BTD in rank-$(L,L,1)$. 

Since the BTD-$(L_r, L_r, 1)$ matches the model for the air dataset, we can use the algorithm BTD-RALS to solve for the air sample tensor to minimize the above function (\ref{eq:costfunction}).

In \cite{Hopke}, the cost function used is
\begin{eqnarray}\label{eq:costfunction1}
\qquad\quad Q=\sum_{i=1}^I \sum_{j=1}^J \sum_{k=1}^K \biggl(\frac{x_{ijk}-\sum_{p=1}^P a_{ijp}b_{kp}}{u_{ijk}} \biggl)^2,
\end{eqnarray}
where $u_{ijk}$ is an uncertainty estimate element in the $i$th species of the $j$th particle size of the $k$th time sample. The procedure of Polissar \cite{Polissar} was used to assign measured data and the associated uncertainties as the input data.


Comparing the functions (\ref{eq:costfunction1}) (element-wise) and (\ref{eq:costfunction}), the cost function in our method does not include the uncertainty estimates. We will use the cost function without the uncertainties $u_{ijk}$ and consider the function (\ref{eq:costfunction1}) as a constraint on our solution. 

According to the paper \cite{Hopke}, there are $9$ sources. Thus, in the BTD model (\ref{eq:costfunction}), we let $P=9$. For the tensor $\mathcal{X} \in \mathbb{R}^{27 \times 3 \times 316}$, the block term decomposition in rank-$(L_r,L_r,1)$ with 9 terms is not essentially unique (see \cite{Lieven2}). Recall that essentially uniqueness indicates the decomposition is unique up to permutation and scaling. Since the tensor data does not satisfy the uniqueness criteria, then we will have multiple solutions from the algorithm. We tested a large number of initial conditions with different $L_p$, $1 \leq p \leq 9$ and found that under the setting of $L_p =9$, $p=1, 2, \dots, 9$, we can find a solution that is consistent with the numerical results found in \cite{Hopke}. 

For the non-negativity constraints on the decomposition, we need to use the BTD-RALS method to solve for the non-negative factor matrices in the problem (\ref{eq:costfunction1}). So we will add non-negativity constraints on the three subproblems. In terms of solving these subproblems, for each least-squares problem with non-negativity constraints, we can use the method in \cite{Lawson}, or we can use the non-negative matrix factorization (NMF) introduced by Lee and Seung \cite{DHS} to solve the each subproblem with constraint.


In this paper, we use the algorithm in Table \ref{tab:btdll1rals} with the method by \cite{Lawson} to obtain the three non-negative factor matrices $\bold{C} $, $\bold{D}$ and $\bold{B}$ in the cost function (\ref{eq:costfunction}). 

\subsection{Result and Discussion}

\begin{figure*}
\begin{center}
\includegraphics[width=.49\textwidth]{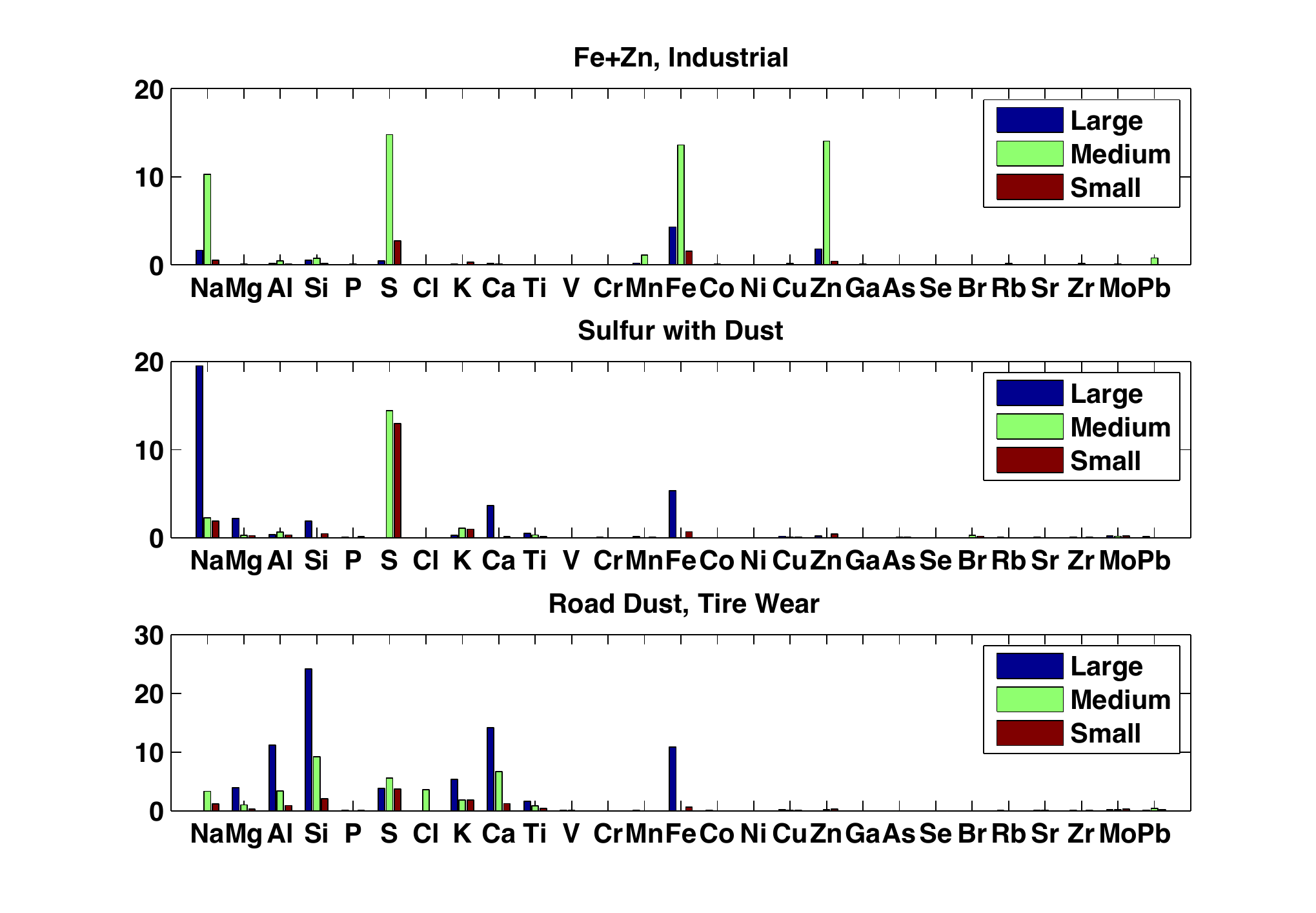}
\hskip2pt
\includegraphics[width=.49\textwidth]{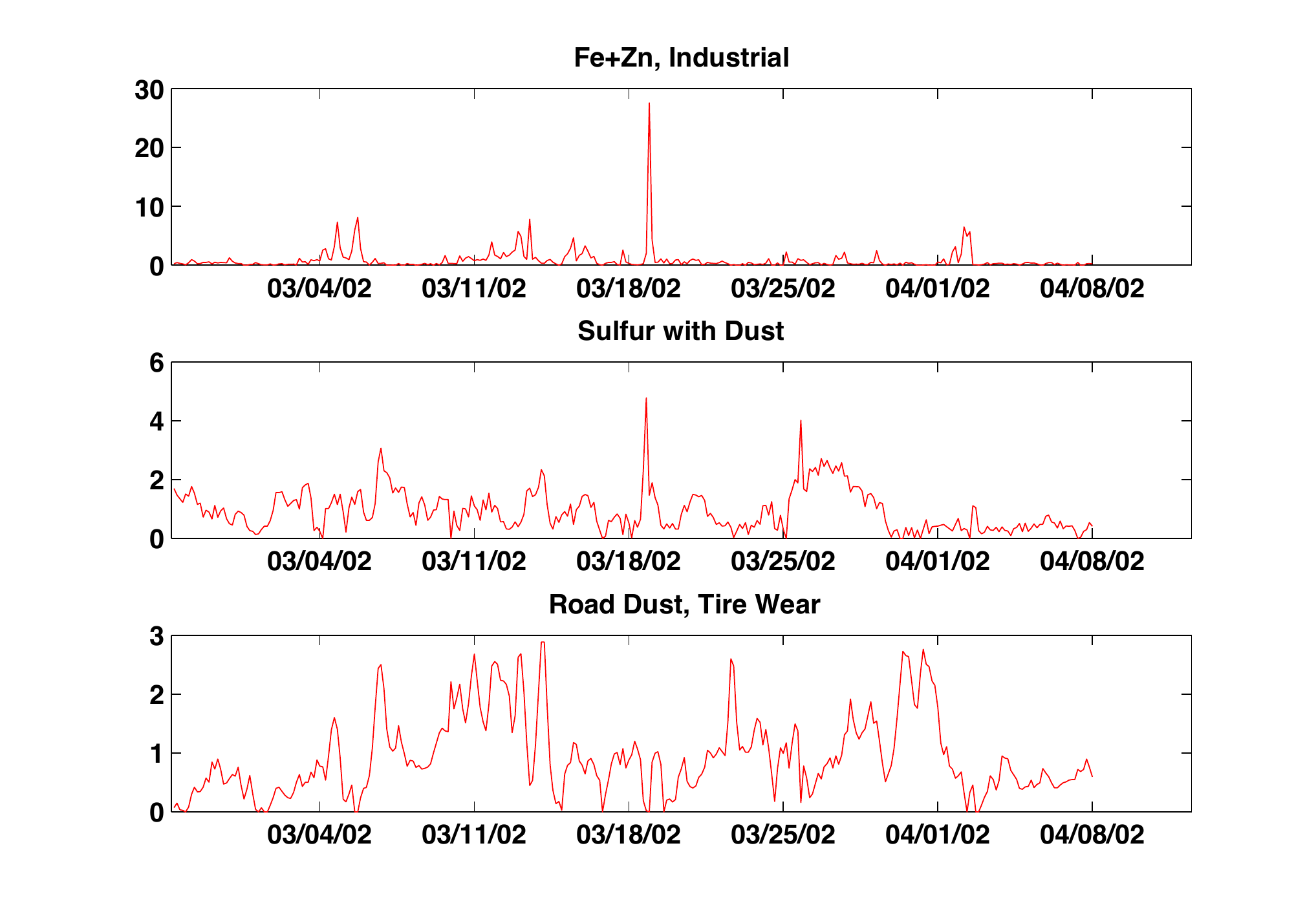}
\vskip2pt
\includegraphics[width=.49\textwidth]{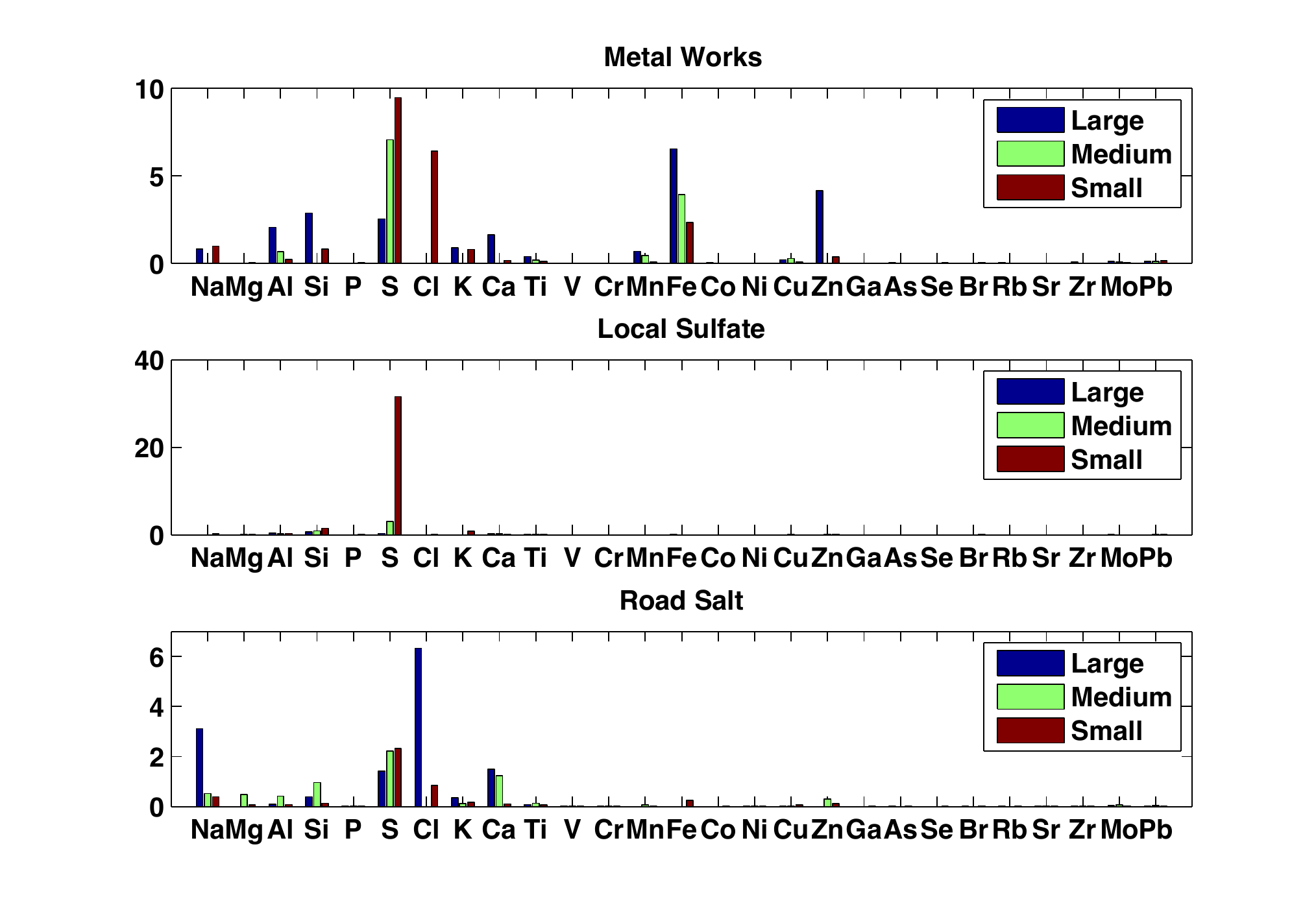}
\hskip2pt
\includegraphics[width=.49\textwidth]{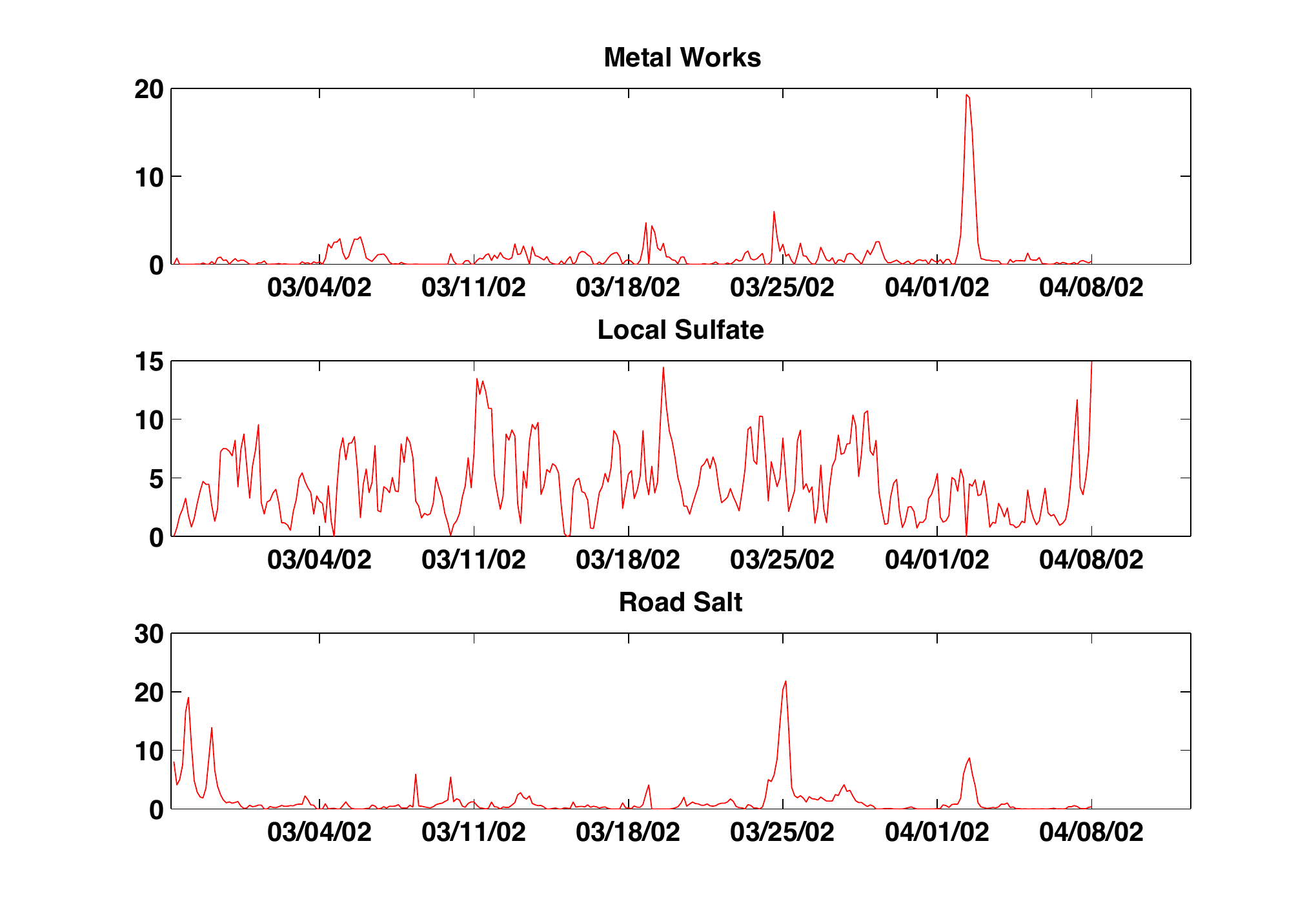}
\vskip2pt
\subfloat[Source profiles for the resolved factors. The $y$-axis denotes the relative elemental concentration. The order of the particle size is Large, Medium, Small from left to right for each element.]{\label{fig:bar}\includegraphics[width=.49\textwidth]{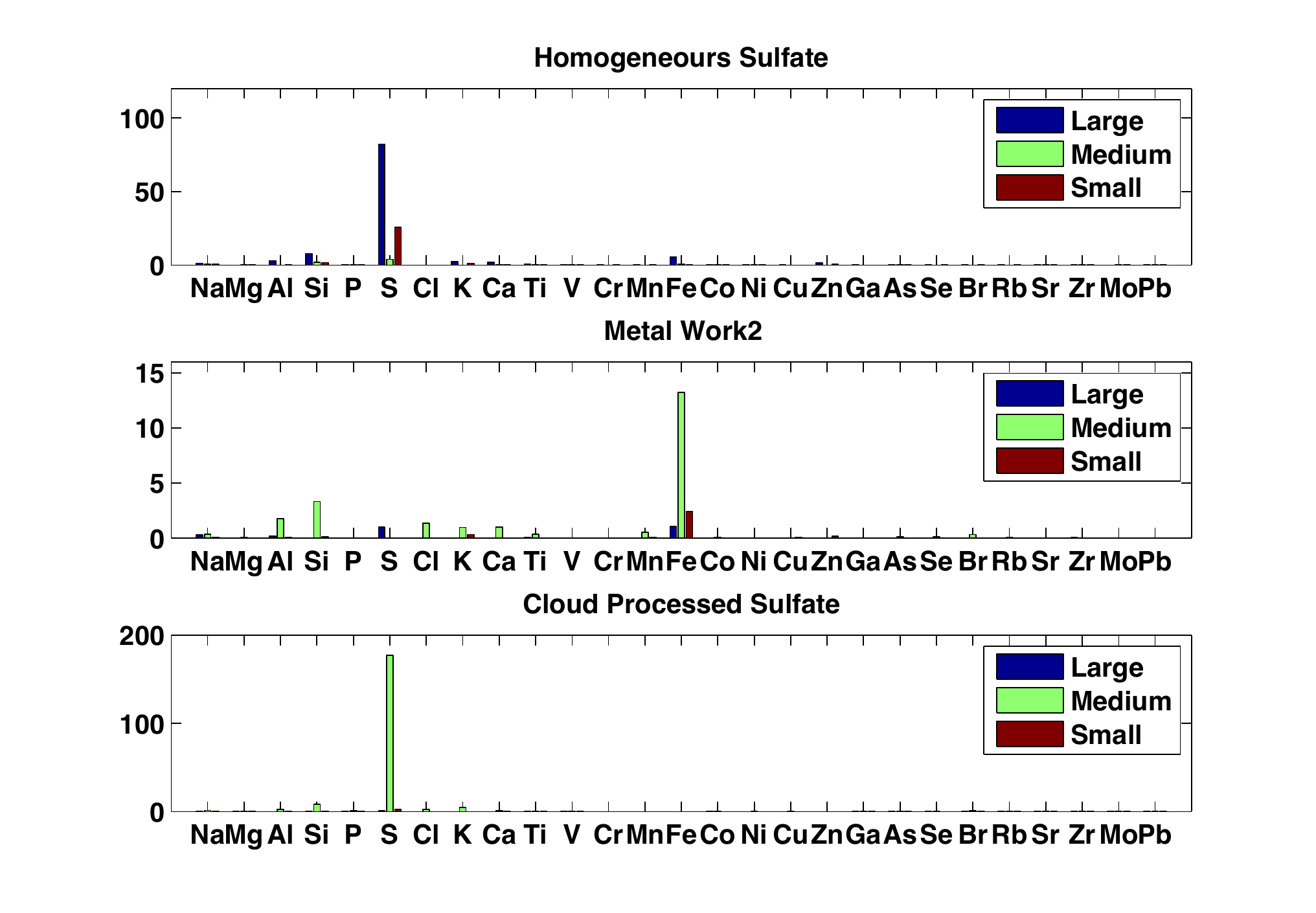}}
\subfloat[The time series of source contributions. ]{\label{fig:time}\includegraphics[width=.49\textwidth]{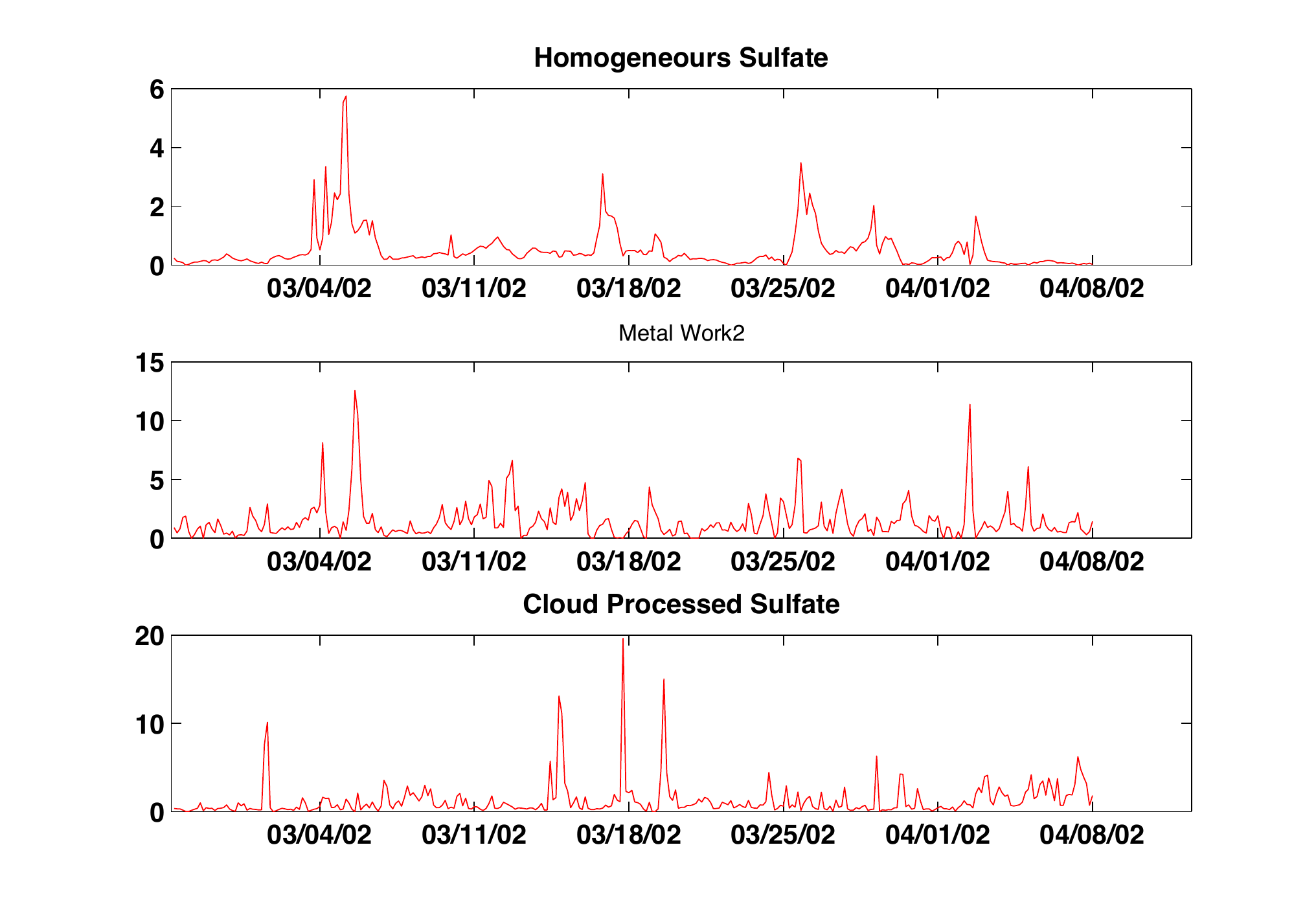}}
\caption{\bf Nine factor sources identified from the air sample.}
\label{eq:result}
\end{center}
\end{figure*}

We apply the BTD-$(9,9,1)$ on the air sample tensor $\mathcal{X}$ (\ref{eq:costfunction}). By using the BTD-RALS algorithm with non-negativity constraints, we can obtain the three factor matrices 

$$\bold{C}=[\bold{C}_1~\cdots ~\bold{C}_9],\;\; \bold{D}=[\bold{D}_1~\cdots~ \bold{D}_9],\;\; \bold{B}=[\bold{b}_1~\cdots~\bold{b}_9],$$
where $\bold{C}_p \in \mathbb{R}^{27 \times 9}$, $\bold{D}_p \in \mathbb{R}^{3 \times 9}$ and $\bold{B} \in \mathbb{R}^{316\times 9}$, $p=1,2,\dots,9$. Therefore, according to the model for the air sample, for each $p$, $\bold{A}_p = \bold{C}_p \cdot \bold{D}_p^T \in \mathbb{R}^{27 \times 3}$ is a matrix and the elements of such matrix are exactly the $a_{ijp}$s (for a fixed $p$) in the equation (\ref{eq:model1}). Furthermore, each $\bold{A}_p$ denotes a source profile. So, for each $\bold{A}_p$ we have a bar plot for one source (see the Figure \ref{fig:bar}). The vector $\bold{b}_p$ in the factor matrix $\bold{B}$ is the vector $b_{kp}$ (for a fixed $p$) in the model equation (\ref{eq:model1}). It denote the contribution of $p$ source in terms of time. Therefore, the Figure \ref{fig:time} expresses the time series of the source contributions.

From the source profile bar plot \ref{fig:bar}, we can figure out the nine different factor sources, they are: Industrial (Fe+Zn), Sulfur with Dust, Road Dust, two types of Metal Works, Local Sulfate, Road Salt, Homogeneously formed Sulfate and Cloud Processed Sulfate. For the explanation for each source, refer to the paper \cite{Hopke}. In the Figure \ref{fig:time}, we can also see the change of each source contribution in time. For example, there are several spikes in the time series of the road source contribution. This indicates additional snowfall or low temperatures where the ice melting. Furthermore, we can also tell the contribution change of other factor sources.

The identification from the bar plots in Figure \ref{fig:bar} is also based the particle sizes. This method provides a more accurate result comparing the classical matrix factorization method. It is seen that in the industrial source, the high concentrations of Fe and Zn are in the middle size rage while the large size rage of Fe and Zn are shown in the metal works.

\begin{figure*}
\begin{center}
\includegraphics[width=.9\textwidth]{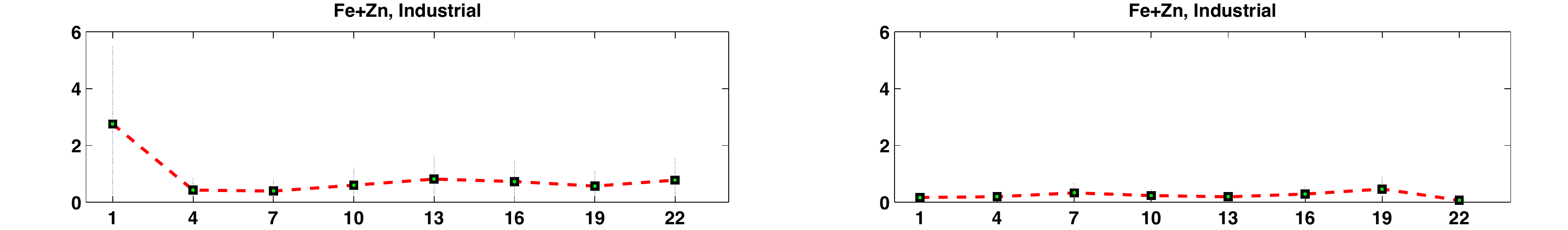}
\includegraphics[width=.9\textwidth]{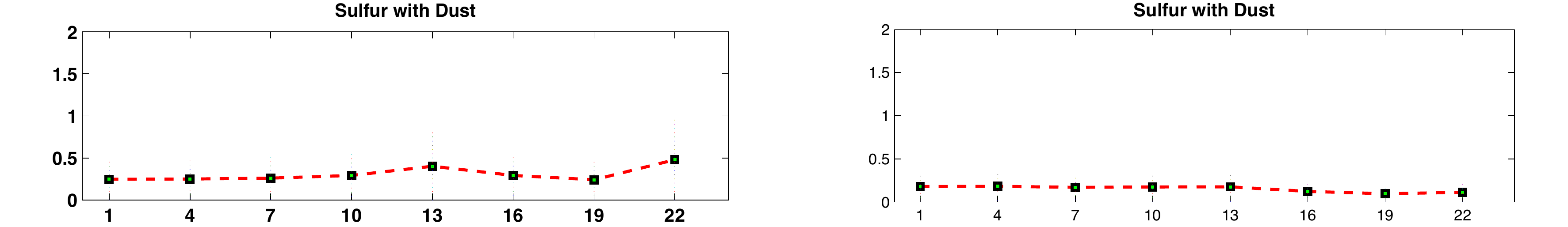}
\includegraphics[width=.9\textwidth]{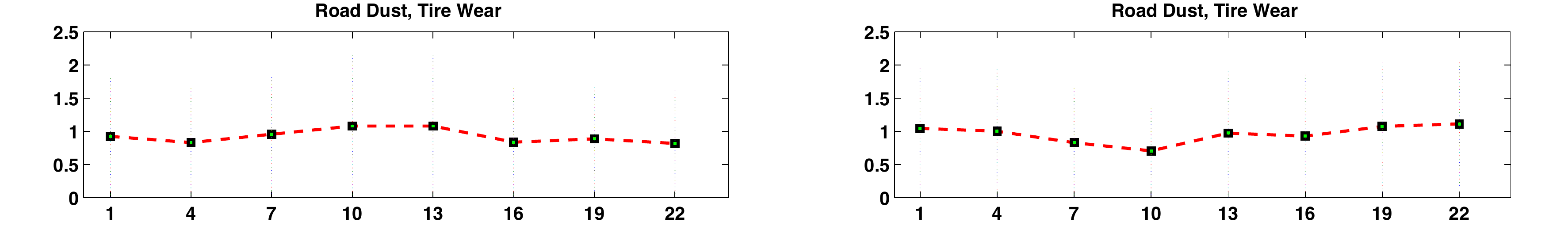}
\includegraphics[width=.9\textwidth]{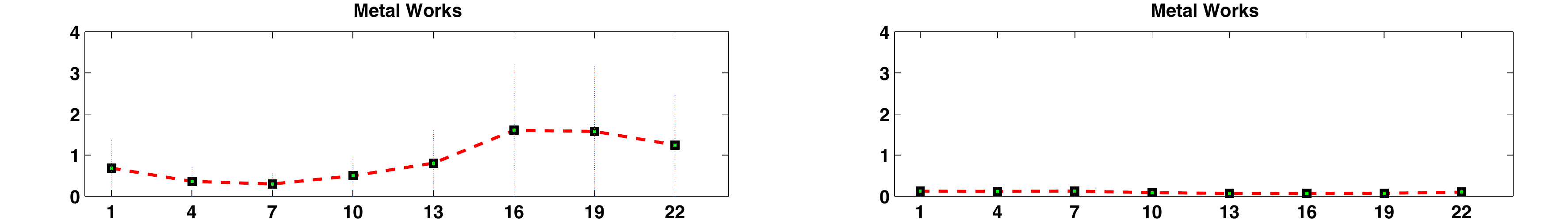}
\includegraphics[width=.9\textwidth]{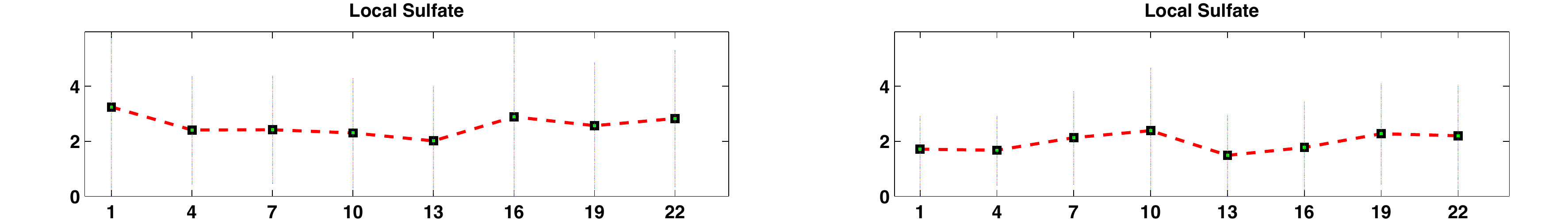}
\includegraphics[width=.9\textwidth]{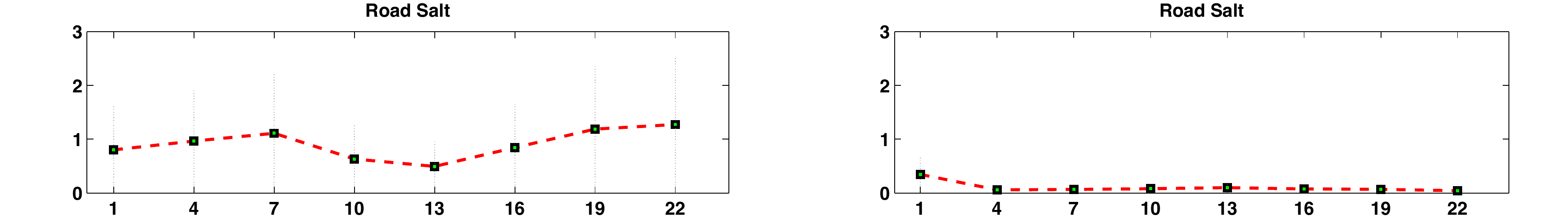}
\includegraphics[width=.9\textwidth]{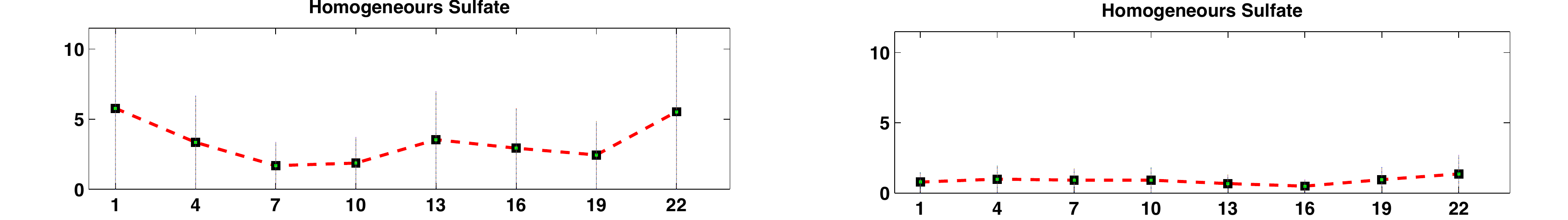}
\includegraphics[width=.9\textwidth]{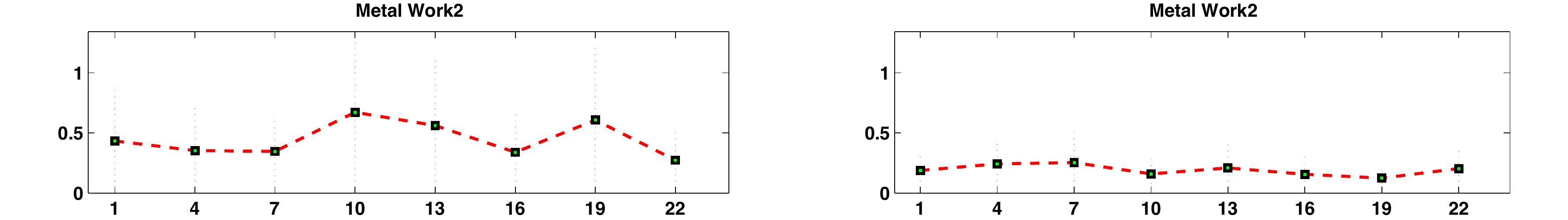}
\includegraphics[width=.9\textwidth]{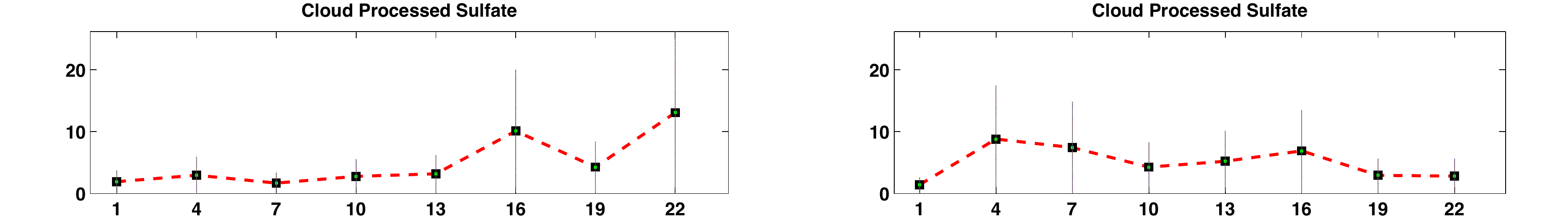}
\caption{\bf Time of a day variations. The left column shows the concentration change in weekdays for each factor source. The right column shows the concentration change in weekends for each factor source}
\label{fig:timeday}
\end{center}
\end{figure*}

We can also analyze the concentration during the weekdays and weekends for each factor source (see the Figure \ref{fig:timeday}). The left nine plots show the concentration change during the weekdays and the corresponding right plots show the concentration change during the weekends for each factor source. In each day, we take the maximum and minimum concentration at the time points 1, 4, 7, 10, 13, 16, 19, 22 and then take the average for each time point. Figure \ref{fig:timeday} shows the factors concentrations of Road Salt, Homogeneous Sulfate and the second type of Metal Work are high during the weekdays and lower during the weekends. This is very reasonable since the weekdays are typically work days when industrial companies are in operation. There is also a longer commuting hours during the weekdays than the weekends, explaining the higher concentration of Road Salt. For some factors like local sulfate, there is no difference between the weekdays and weekends which indicates the people's activities are not a major effect on these environmental factors. 

With the Figures \ref{fig:bar}, \ref{fig:time} and \ref{fig:timeday} and above discussion, we show that the BTD-RALS method can successfully identify the nine different factor sources and also provides the time contribution of each factor. 


\subsection{BTD-RALS on the noisy data}
Here we give a numerical experiment to show the BTD-RALS method works on the third-order tensor data with noise. We generate tensors $\widetilde{\mathcal{T}} \in \mathbb{R}^{5 \times 6 \times 7}$ in the following way:
\begin{eqnarray}\label{eq:noise}
\widetilde{\mathcal{T}} = \frac{\mathcal{T}}{\Vert \mathcal{T} \Vert} + \sigma_N \frac{\mathcal{N}}{\Vert \mathcal{N}\Vert},
\end{eqnarray}
where $\mathcal{T}$ has the block term decomposition in rank-$(2,2,1)$ with $R=3$ in the equation (\ref{eq:btd2}) so that $\bold{A}_r \in \mathbb{R}^{5 \times 2}$, $\bold{B}_r \in \mathbb{R}^{6 \times 2}$ and $\bold{c}_r \in \mathbb{R}^{7 \times 1}$, $r=1, 2, 3$. The second term in (\ref{eq:noise}) is the noise term and the parameter $\sigma_N$ controls the noise level. The entries of $\bold{A}=[\bold{A}_1~\bold{A}_2~\bold{A}_3]$, $\bold{B}=[\bold{B}_1~\bold{B}_2~\bold{B}_3]$, $\bold{C}=[\bold{c}_1~\bold{c}_2~\bold{c}_3]$ and the tensor $\mathcal{N}$ are drawn from a zero-mean unit-variance Gaussian distribution. 

By using BTD-RALS method, a Monte Carlo experiment of BTD-$(2,2,1)$ with $R=3$ on $\widetilde{\mathcal{T}}$ consisting 50 runs is carried out. The algorithm is initialized with three random starting values.

\begin{figure}[h!]
\begin{center}
\includegraphics[width=.45\textwidth]{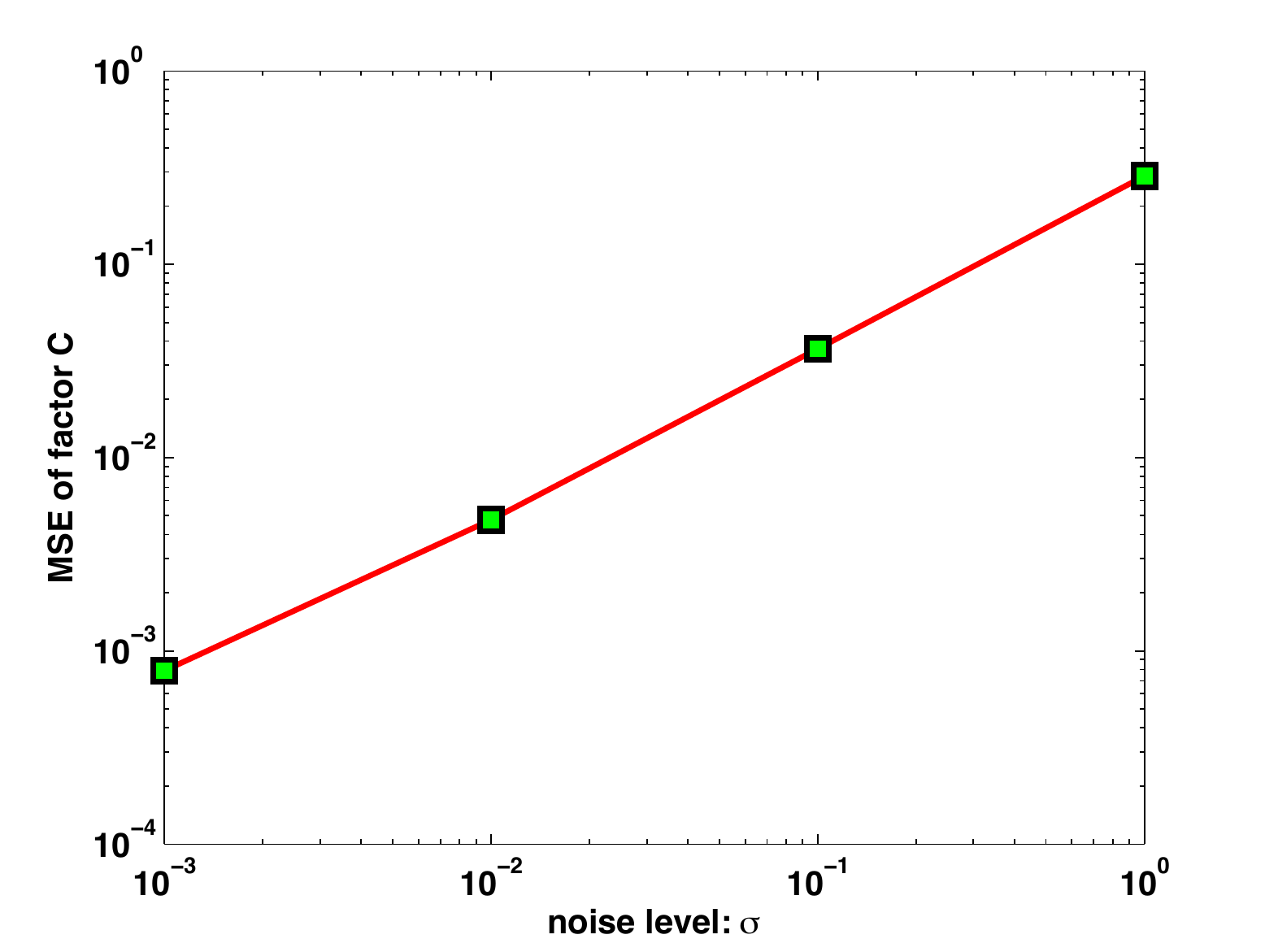}
\caption{Median relative error v.s. noise level}
\label{fig:noise}
\end{center}
\end{figure}

We measure the relative error $e=\Vert \bold{C}-\widetilde{\bold{C}} \Vert/ \Vert \bold{C}\Vert$, where $\widetilde{\bold{C}}$ is the approximation of $\bold{C}$, optimally permuted and scaled. The median results are plotted in Figure \ref{fig:noise}. 
It is seen that with low noise levels, average error in $\bold{C}$ increases proportionally to noise level. 

\section{Conclusion}
The method BTD-RALS is presented in the application of identifying the factor sources of the collected air sample. The dataset is formed into a third order tensor and the data model is written into a block term decomposition in rank-$(L,L,1)$. We apply a regularized alternating method to solve the block term decomposition and obtain the resulting factor matrices providing the source profiles and source contributions correctly. In addition, we show that regularized method is efficient than the classical alternating method numerically and can converge to a stationary point of the original BTD cost function under some assumption. The BTD-RALS algorithm is also tested on the random data with different noise levels, which shows the average error in the third factor matrices $\bold{C}$ increases proportionally to noise level.


\end{document}